ON TWO COMBINATORIAL INEQUALITIES THAT EXPLAIN THE BLIMPY SHAPE OF HEADY-*S* AND TAILY-*S* BIT STRINGS

**Bruce Levin** BL6@COLUMBIA.EDU
*Department of Biostatistics*
*Mailman School of Public Health*
*Columbia University*
*New York, NY 10032, USA*

**Abstract**

We prove two combinatorial inequalities introduced in our prior study of the "blimpy" graphical shape of the number of bit strings with a given score under an interesting scoring system. Generating functions are used to establish the inequalities, which in turn imply two of the salient graphical features, unimodality near the zero score and shape asymmetry for positive versus negative scores. One inequality provides a lower bound on the expected value of a discrete random variable with probabilities proportional to a product of two binomial coefficients. The other inequality states that the expected value with respect to peri-central binomial coefficients of other binomial coefficients lying on an oblique ray in Pascal's triangle exceeds the expected value along an adjacent parallel ray to its left.

**Keywords**

Combinatorial inequalities; discrete convexity; expectation inequalities; generating functions.

**Subject Classifications:** 05A20; 60C05.

## 1. INTRODUCTION

Levin (2025, hereinafter "L25") studied the "blimpy" shape of the graph of the logarithm of the number of bit strings of length *n* which possess certain "score" values, when plotted against those values, under an interesting scoring system for bit strings. The motivation for the scoring system begins with the following probability puzzle posed by Daniel Litt of the University of Toronto.

> Alice and Bob flip a (fair) coin 100 times. Anytime there are two heads in a row, Alice gets a point; when a head is followed by a tail, Bob gets a point. So in the sequence THHHT, Alice gets two points and Bob gets one. Who is more likely to win?
>
> [Levin (2024, hereinafter L24), paraphrasing Litt (2024) and Klarreich (2024)]

The score in question is the net difference between Alice's and Bob's point totals after *n* tosses. The puzzle's answer is counterintuitive—anytime the coin is tossed more than twice, Bob is more likely to win than Alice. Attempting to understand why this is so at an intuitive level leads to a quandary of heuristic reasoning. Most people begin with an appeal to symmetry by reasoning that the coin is fair, so after any toss of a head, both players have an equal chance of scoring a point on the next toss, leading to the wrong conclusion that the odds on winning should be the same for Alice and Bob. Some people argue heuristically that Alice can accrue a given number of points with *overlapping* (1,1) bit pairs in fewer tosses than Bob can accrue with only *non-overlapping* (1,0) bit pairs. For example, Alice can make three points in four tosses at the start of the game, but it takes Bob at least six tosses to make three points. That seems to give Alice an advantage, but this reasoning also leads to a wrong conclusion. A countervailing heuristic is that it takes Bob two fewer tosses on average to reach his first point than it takes Alice to reach hers (four versus six), advantage Bob. Of course, reaching different conclusions when comparing event probabilities of random variables versus their averages is a hallmark of asymmetrical distributions. Intuition aside, it is easy to see by simple enumeration that no matter what the score may be after $n$–3 tosses, Bob has more ways to *improve* his score by game's end than Alice does. How to turn this observation into a rigorous proof that Bob does have the advantage is not at all obvious. L24 does this with an elementary inductive argument (with no generating functions or asymptotic calculations). Section 1 of L25 further discusses the asymmetric odds on winning as have several other authors, including Ekhad and Zeilberger (2024), Segert (2024), Basdevant et al. (2024), and Janson et al. (2025). Of interest to readers of this journal, a simple sequential version of the game is perfectly symmetrical and fair (see L25, Section 1), though a sequential stopping time formulation to address the original asymmetrical fixed-sample size game runs into difficulties, mentioned below, due to the need to condition on certain required events (see L25, Section 4). Other interesting asymmetries arise, though, which is what motivated L25 and the present work.

We call a binary string $x^{(n)} = (x_1, ..., x_n) \in \{0,1\}^n$ with $x_n = 1$ a *heady* bit string (of length $n$) or a *taily* bit string if $x_n = 0$. The *score* of $x^{(n)}$ is then $S(x^{(n)}) = \sum_{i=1}^{n-1} I[x_i = x_{i+1} = 1] - \sum_{i=1}^{n-1} I[x_i = 1,\ x_{i+1} = 0]$. For a given integer $s$, we call a heady (respectively, taily) bit string with score $S(x^{(n)}) = s$ a *heady-s* (respectively, *taily-s*) bit string. Figure 1.1 displays the exact number of heady-$s$ bit strings (left panel) and taily-$s$ bit strings (right panel) as a function of $s$ (center column) for the case $n=100$. The shape of the figure exhibits *unimodality*, with modal value at $s=0$ for heady strings and at $s=-1$ for taily strings (see the enlargement), and *asymmetry* in $s$, with elongation for positive $s$ and foreshortening for negative $s$. Paraphrasing L25, people tend to see different shapes in the image, like a leaf or an egg or a teardrop, but the diagram reminds us of a dirigible. What explains that "blimpy" shape?

[FIGURE 1 HERE]

Let $H_s(n)$ denote the number of heady-$s$ bit strings of length $n$. Then we have

$$H_s(n) = \sum_{k=0\vee(-s)}^{\lfloor (n-s-1)/3 \rfloor} \binom{2k+s}{k} \binom{n-s-1-2k}{k}, \tag{1.1}$$

which is positive for $n \geq 2$ and $s$ in the range $-\lfloor (n-1)/2 \rfloor \leq s \leq n-1$ and zero elsewhere (L24, L25). The following theorem specifies the unimodal shape of $H_s(n)$ (restating L25, Theorem 1.2).

Theorem 1.1 (unimodality of $H_s(n)$ in $s$ for given $n$).

(a) *For any* $s \ge 0$, $H_s(n) \ge H_{s+1}(n)$. (1.2)

(b) *For any* $s \le 0$, $H_s(n) \ge H_{s-1}(n)$. (1.3)

The number of taily-$s$ bit strings, $T_s(n)$, has an analogous formula and theorem, but an elementary argument shows that Theorem 1.1 actually implies the corresponding unimodality of $T_s(n)$ at $s = -1$ (L25, Lemma 1.1). For this reason, going forward we focus exclusively on the function $H_s(n)$.

The not-quite-complete proof we gave of Theorem 1.1 did not use generating functions, in part to reveal special properties of $H_s(n)$ and related discrete distributions which might already imply unimodality in $s$, properties such as functional convexity, monotonic expectations of stopping times, and pairwise point asymmetry of certain discrete functions. Those methods succeeded up to a point, but left some gaps which we bridged by stating additional properties conjectured but not proven to hold. *Our purpose now is to prove Theorem 1.1 completely with generating functions.*

The asymmetrical shape of Figure 1.1 with respect to $s$ emerges from the truncation of $k$ to $\max(0,-s)$ in summation (1.1). It also manifests in how the necessary and sufficient conditions for unimodality differ for non-negative versus non-positive $s$, and in the two combinatorial inequalities that we must marshal to meet those conditions. The two key inequalities appear below at (1.8) and (1.18). We turn to these topics next, appreciating in passing how nicely large-scale structure emerges from local combinatorial properties of binary sequences and related systems.

For the convenience of the reader, Appendix C provides a glossary of the notation used throughout the article. In the case of non-negative $s$, the following condition is necessary and sufficient for Theorem 1.1(a) to hold (restating L25, Lemma 2.1).

Lemma 1.1. *For any $n$ and $0 \le s < n-1$, let $n_s = n-s-1$ and $\kappa = \kappa_s = \lfloor n_s/3 \rfloor$. To avoid trivialities, assume that $n_s \ge 3$ so that $\kappa \ge 1$. Then $H_s(n) \ge H_{s+1}(n)$ if and only if*

$$2E^{(1)}_{s,n}\left[\binom{n_s-1-2K}{K-1}\right] + \Delta^{(1)}_{s,n} \ge E^{(1)}_{s,n}\left[\binom{n_s-2K}{K}\right] \ge 2E^{(1)}_{s,n}\left[\binom{n_s-1-2K}{K}\right] - \Delta^{(1)}_{s,n}, \quad (1.4)$$

*where $K$ is a discrete random variable with probabilities proportional to the leading binomial coefficient in $H_s(n)$, namely,*

$$P^{(1)}_{s,n}[K=k] = \binom{2k+s}{k} \Big/ \sum_{j=1}^{\kappa} \binom{2j+s}{j} \quad \text{for } k = 1,\ldots,\kappa; \quad (1.5)$$

*where the expectations are taken with respect to (1.5); and where*

$$\Delta_{s,n}^{(1)} = E_{s,n}^{(1)}\left[\binom{n_s-1-2K}{K}\left(\frac{s+1}{K+s+1}\right)\right]. \quad \square \tag{1.6}$$

It is not obvious how to establish that (1.4) holds working under (1.5), but a change of measure to absorb the second binomial coefficient in with the first does allow us to show alternatively that $H_s(n) \ge H_{s+1}(n)$ *if and only if* $E_{s,n}^{(2)}[f(K)] \ge 0$, where now the expectation is with respect to

$$P_{s,n}^{(2)}[K=k] \propto \binom{2k+s}{k}\binom{n_s-2k}{k} \text{ for } k=1,\ldots,\kappa, \tag{1.7}$$

and where $f(x) = \left(\frac{x}{x+s+1}\right)\left(\frac{4x+s+1-n_s}{n_s-2x}\right)$ for $0 \le x \le n_s/3$. For larger values of $s$ such that $n_s \le s+1$, $f(x)$ is non-negative over the entire interval, in which case we are done. For smaller values of $s$ such that $n_s > s+1$, $f(x)$ has a zero at $x^* = \{n_s-(s+1)\}/4$, is increasing there, and is convex over the entire interval (L25, Section 3). Thus it would suffice by Jensen's inequality if

$$E_{s,n}^{(2)}[K] \ge \frac{n_s-(s+1)}{4} = \frac{n-2(s+1)}{4} = x^*, \tag{1.8}$$

for then $E_{s,n}^{(2)}[f(K)] \ge f(E_{s,n}^{(2)}[K]) \ge f(x^*) = 0$ since $f$ increases to the right of $x^*$ by convexity.

*We prove inequality (1.8) in Theorem 2.1 below.* In the present case of non-negative $s$, the somewhat opaque necessary and sufficient condition (1.4) or its equivalent $E_{s,n}^{(2)}[f(K)] \ge 0$ shepherds us into first showing the stronger inequality $E_{s,n}^{(2)}[K] \ge \{n_s-(s+1)\}/4$, from which (1.4) then follows. The reverse logical direction arises in the case of non-positive $s$, which we consider next.

For $s \le 0$, in what follows we write $u = -s \ge 0$ to minimize the number of minus signs and will use the shorthand notation $n_u = n+u = n-s$ (rather than $n_s = n-s-1 = n+u-1$), though we continue to use subscript $s$ in $H_s(n)$ over the entire range of $s$ for consistency. Now (1.1) can be written

$$H_s(n) = \binom{n-u-1}{u} + \sum_{k=u+1}^{\lfloor (n_u-1)/3 \rfloor} \binom{2k-u}{k}\binom{n_u-1-2k}{k} = \binom{n-u-1}{u} + \sum_{k=u+1}^{\lfloor (n_u-1)/3 \rfloor} a_k(n,u), \text{ say,} \tag{1.9}$$

and

$$H_{s-1}(n) = \sum_{k=u+1}^{\lfloor n_u/3 \rfloor} \binom{2k-u-1}{k}\binom{n_u-2k}{k} = \sum_{k=u+1}^{\lfloor n_u/3 \rfloor} b_k(n,u), \text{ say,}$$

and we want to show $H_s(n) \ge H_{s-1}(n)$. Under the conditions

$$n \ge 7 \text{ and } 0 \le u < \lfloor (n-1)/2 \rfloor - 1 = \lfloor (n-3)/2 \rfloor, \tag{1.10}$$

the leading term in (1.9) turns out to be superfluous, i.e., it is already the case that

$$\sum_{k=u+1}^{\lfloor n_u/3 \rfloor} a_k(n,u) \geq \sum_{k=u+1}^{\lfloor n_u/3 \rfloor} b_k(n,u)\,, \tag{1.11}$$

which of course implies (1.3). We may use the same upper limit of summation on both sides of (1.11) because when $\lfloor (n_u-1)/3 \rfloor = \lfloor n_u/3 \rfloor - 1$, i.e., when $n \equiv 0 \bmod 3$, the last addend on the left of (1.11) equals zero. For values of $n$ and $u$ not satisfying (1.10), elementary arguments directly show (1.3) holds, so we assume (1.10) henceforth.

The necessary and sufficient condition for (1.11) is now as follows (L25, Lemma 2.5).

Lemma 1.2. *For given n and u satisfying (1.10), inequality (1.11) holds if and only if*

$$2E_{u,n}^{(3)}\left[\binom{n_u-1-2K}{K}\right] + \Delta_{u,n}^{(3)} \;\geq\; E_{u,n}^{(3)}\left[\binom{n_u-2K}{K}\right] \;\geq\; 2E_{u,n}^{(3)}\left[\binom{n_u-1-2K}{K-1}\right] - \Delta_{u,n}^{(3)}\,, \tag{1.12}$$

*where the expectations are taken with respect to*

$$P_{u,n}^{(3)}[K=k] = \binom{2k-u-1}{k} \Big/ \sum_{j=1}^{\kappa} \binom{2j-u-1}{j} \quad \textit{for } k = u+1,\ldots,\kappa_u = \lfloor n_u/3 \rfloor, \tag{1.13}$$

*and*

$$\Delta_{u,n}^{(3)} = E_{u,n}^{(3)}\left[\binom{n_u-1-2K}{K}\left(\frac{u}{K-u}\right)\right]. \quad \square \tag{1.14}$$

Another change of measure to include the second binomial coefficient, plus a change of variables from $K$ to $K' = K-u$ over the range $k' = 1,\ldots,\kappa_u' = \kappa_u - u = \lfloor (n-2u)/3 \rfloor$, allows us to show that (1.12) *is fully equivalent* to the following inequality which, after dropping the primes, we can write as

$$E_{u,n}^{(4)}\left[K^{-1}\right] \geq \frac{4}{n-2u}\,, \tag{1.15}$$

where now for $k = 1,\ldots,\kappa = \lfloor (n-2u)/3 \rfloor$,

$$P_{u,n}^{(4)}[K=k] \propto \binom{2k+u-1}{k-1}\binom{n-u-1-2k}{k+u-1}. \tag{1.16}$$

As $K^{-1}$ is convex, Jensen's inequality would further imply the sufficient condition

$$E_{u,n}^{(4)}[K] \leq (n-2u)/4\,, \tag{1.17}$$

but now a problem arises, because (1.17) is definitely *not true* for $u$=0, and, although (1.15) does hold for all $u \geq 0$ as we shall see, the methods we used in L25 to establish (1.15), or (1.17) for $u > 0$, gave a complete proof for all *n only* when $u > 2$, leaving gaps for certain subsets of *n* when $u$=0, 1, and 2. However, we need not first show (1.15) for such *u*. Instead, we can show that for *n* sufficiently large,

$$E_{u,n}^{(3)}\left[\binom{n_u-1-2K}{K}\right] \geq E_{u,n}^{(3)}\left[\binom{n_u-1-2K}{K-1}\right]. \tag{1.18}$$

Inequality (1.18) implies condition (1.12), because *a fortiori* from (1.18) we have

$$2E_{u,n}^{(3)}\left[\binom{n_u-1-2K}{K}\right] + \Delta_{u,n}^{(3)} \geq 2E_{u,n}^{(3)}\left[\binom{n_u-1-2K}{K-1}\right] - \Delta_{u,n}^{(3)} \tag{1.19}$$

since $\Delta_{u,n}^{(3)}$ in (1.14) is non-negative, and the middle term of (1.12) lies between the two sides of (1.19) because the former is a simple average of the latter by Pascal's identity. Condition (1.15) then follows from its equivalence to (1.12). *We prove inequality (1.18) under condition (1.10) when u=0, 1, or 2 and n is sufficiently large using generating functions in Theorem 3.1 below.*

Note that we couldn't argue in this direction for non-negative *s* because it is *not generally true that* $E_{s,n}^{(1)}\left[\binom{n_s-1-2K}{K-1}\right] \geq E_{s,n}^{(1)}\left[\binom{n_s-1-2K}{K}\right]$, so the above *a fortiori* argument wouldn't work. Indeed, the always-positive $\Delta_{s,n}^{(1)}$ in (1.6) is essential for (1.4) to hold.

Inequality (1.18) has some combinatorial interest of its own. In the case of $u=0$, for example, (1.18) states the subtle fact that for $n \geq 7$, the weighted average of $\binom{n-1-2k}{k}$ over $k=1,\ldots,\lfloor n/3 \rfloor$ with weights proportional to $\binom{2k-1}{k}$ is always at least that of $\binom{n-1-2k}{k-1}$. In general, the averaging takes place over parallel oblique lines through points in Pascal's triangle.

## 2. THE CASE OF NON-NEGATIVE SCORES

In this section we assume $s \geq 0$. It will be convenient to use the index *n* to refer to what we previously called $n_s = n-s-1$ rather than the original bit string length *n*. Then inequality (1.8) becomes

$$E_{s,n}[K] \geq (n-s-1)/4 \tag{2.1}$$

with respect to the probability function

$$P_{s,n}[K=k] \propto \binom{2k+s}{k}\binom{n-2k}{k} \quad \text{for } k=1,\ldots,\kappa=\lfloor n/3 \rfloor. \tag{2.2}$$

The expectation in (2.1) is what we would have written as $E_{s,n+s+1}^{(2)}[K]$ in Section 1. Also, instead of the constraint $s<n-1$ as in Lemma 1.1, we may now simply let *s* be any non-negative integer.

<u>Theorem 2.1</u>. *For* $n \geq 3$, $E_{s,n}[K] \geq (n-s-1)/4$. *The inequality is strict for all such n when* $s>0$, *and when s=0, it is strict for all such n except n=5 and 8, where* $E_{s,n}[K]=(n-s-1)/4$.

Proof. If $3 \le n < 6$, then $\kappa = 1$ and (2.2) degenerates to $K=1$, in which case (2.1) holds since $(n-s-1)/4 \le 1$ for any $s \ge 0$. In the degenerate case $n=5$ with $s=0$, $E_{0,5}[K]=(n-1)/4=1$. So, without loss of generality, we may now assume $n \ge 6$. In the non-degenerate case $n=8$ with $s=0$, $E_{0,8}[K]=(1\cdot 2\cdot 6+2\cdot 6\cdot 6)/(2\cdot 6+6\cdot 6)=84/48=7/4=(n-1)/4$. Next we show strict inequality holds for all other $s$ and $n \ge 6$.

First consider the case $s=0$ and for any non-negative $n$, define the sums

$$A_n = \sum_{k\ge 0}\binom{2k}{k}\binom{n-2k}{k} \tag{2.3}$$

and

$$B_n = \sum_{k\ge 0}k\binom{2k}{k}\binom{n-2k}{k}. \tag{2.4}$$

The sums are finite since all terms vanish for $k > n/3$ under the usual convention for binomial coefficients [namely, that they are zero if the lower integer index is negative or if it is greater than the upper index when that is an integer; in all other cases with real $x$, $\binom{x}{i} = x^{(i)}/i! = x(x-1)\cdots(x-i+1)/i!$].

Then from (2.2) with $s=0$, $E_{0,n}[K] = B_n/(A_n-1) > B_n/A_n$ and we will show that

$$B_n/A_n > (n-1)/4 \tag{2.5}$$

for all $n \ge 3$ except for $n=5$ and 8, where the inequality is reversed [though as mentioned above, equality in (2.1) still holds].

We derive the generating functions for sequences $A_n$ and $B_n$ in Appendix A and there we apply the transfer theorems of singularity analysis (Flajolet and Sedgewick, 2009) to show that for large $n$,

$$A_n = \frac{2^n}{\sqrt{\pi n}}\{1 - \frac{1}{16n} + \frac{121}{512n^2} + O(n^{-3})\} \tag{2.6}$$

and

$$B_n = \frac{2^n\sqrt{n}}{4\sqrt{\pi}}\{1 - \frac{13}{16n} - \frac{159}{512n^2} + O(n^{-3})\}. \tag{2.7}$$

Thus,

$$\frac{B_n}{A_n} = \frac{n}{4}\left\{\frac{1-\frac{13}{16n}-\frac{159}{512n^2}+O(n^{-3})}{1-\frac{1}{16n}+\frac{121}{512n^2}+O(n^{-3})}\right\} = \frac{n}{4}\left\{1-\frac{3}{4n}+O(n^{-2})\right\} = \frac{n}{4}-\frac{3}{16}+O(n^{-1}) = \frac{n-1}{4}+\frac{1}{16}+O(n^{-1}), \tag{2.8}$$

so that $(B_n/A_n)-(n-1)/4$ exceeds any positive value less than 1/16 for sufficiently large $n$. The order $1/n$ term on the right-hand side of (2.8) is small ($-\frac{19}{128n}$; see Appendix A) and suggests that $B_n/A_n$ could exceed $(n-1)/4$ even for small *n*. Indeed, exact numerical evaluation shows that $B_n / A_n > (n-1)/4$ for all *n*=4, 6, 7, and $9 \le n \le 500$, extending well into the asymptotic regime.

Remark 2.1. All "exact" numerical evaluations in this paper were carried out in the APL+WIN Version 19.0.01 computing environment (APLNow, LLC, Brielle, NJ) in double-precision arithmetic except for the entries of Figure 1, which were prepared using arbitrary-precision arithmetic, so are entirely exact. APL programs for the calculations are available from the author upon request.

For $s>0$, *we claim that* $E_{s,n}[K] > E_{s-1,n}[K]$ *for all non-degenerate* $K$ *with* $n \ge 6$. Let $g(k) = \binom{2k+s}{k} \Big/ \binom{2k+s-1}{k}$. Then

$$E_{s,n}[K] = \frac{\sum_{k\ge1} k\binom{2k+s}{k}\binom{n-2k}{k}}{\sum_{k\ge1}\binom{2k+s}{k}\binom{n-2k}{k}} = \frac{\sum_{k\ge1} kg(k)\binom{2k+s-1}{k}\binom{n-2k}{k}}{\sum_{k\ge1} g(k)\binom{2k+s-1}{k}\binom{n-2k}{k}} = \frac{E_{s-1,n}[Kg(K)]}{E_{s-1,n}[g(K)]} > E_{s-1,n}[K]$$

if and only if the covariance $Cov_{s-1,n}(K, g(K))$ between *K* and $g(K)$ is positive. But

$$g(k) = \frac{(2k+s)!k!(k+s-1)!}{(k+s)!k!(2k+s-1)!} = \frac{2k+s}{k+s} = 1 + \frac{k}{k+s}$$

is an increasing, non-constant function of *k* (because $\kappa > 1$), so by Chebychev's covariance inequality, $Cov_{s-1,n}(K, g(K)) > 0$, whence $E_{s,n}[K] > E_{s-1,n}[K]$ as claimed. It then follows by induction on *s* that $E_{s,n}[K] > E_{s-1,n}[K] \ge (n-s)/4 > (n-s-1)/4$. Moreover, $E_{s,n}[K] \ge (n-1)/4$ for all $s \ge 0$.

This concludes the proof of Theorem 2.1. □

## 3. THE CASE OF NON-POSITIVE SCORES

In this section we assume $u = -s \ge 0$ and we revert back to the original meaning of *n* as the bit string length. We assume condition (1.10), namely $n \ge 7$ and $0 \le u < \lfloor (n-3)/2 \rfloor$, and work with the $P^{(3)}$ probability function of $K' = K - u$ taking values $k' = 1, \ldots \kappa' = \kappa_u - u$. Dropping the primes, this is now

$$P_{u,n}[K=k] \propto \binom{2k+u-1}{k+u} = \binom{2k+u-1}{k-1} \text{ for } k = 1, \ldots, \kappa = \lfloor (n-2u)/3 \rfloor. \tag{3.1}$$

Then (1.18) becomes

$$E_{u,n}\left[\binom{n-u-1-2K}{K+u}\right] \ge E_{u,n}\left[\binom{n-u-1-2K}{K+u-1}\right]. \tag{3.2}$$

Define the sums

$$C_{u,n} = \sum_{k\ge1} \binom{2k+u-1}{k-1}\binom{n-u-1-2k}{k+u} \tag{3.3}$$

and

$$D_{u,n} = \sum_{k\ge1} \binom{2k+u-1}{k-1}\binom{n-u-1-2k}{k+u-1}. \tag{3.4}$$

As in the previous section, the sums are finite since all terms vanish for $k > (n-2u)/3$. Then (3.2) holds if and only if $C_{u,n} \ge D_{u,n}$, which we next show for $u$=0, 1, and 2. Only these cases are required since we previously showed that sufficient condition (1.17) holds for all $u > 2$ (L25, Lemmas 6.1-6.8).

The generating function (g.f.) $G_u^{(D)}(x) = \sum_{n\ge3} D_{u,n}x^n$ for the $D_{u,n}$ sequence appears in Appendix B and there we show that the g.f. for the $C_{u,n}$ sequence is related to $G_u^{(D)}(x)$ by $G_u^{(C)}(x) = \frac{x}{1-x}G_u^{(D)}(x)$. It follows that the generating function of $C_{u,n} - D_{u,n}$ is

$$\begin{aligned}\left(\frac{x}{1-x}-1\right)G_u^{(D)}(x) &= (-1+x+x^2+\cdots)(D_{u,3}x^3 + D_{u,4}x^4 + D_{u,5}x^5+\cdots) \\ &= -D_{u,3}x^3 + (D_{u,3}-D_{u,4})x^4 + \cdots + \left\{\left(\sum_{j=3}^{n-1} D_{u,j}\right) - D_{u,n}\right\}x^n + \cdots .\end{aligned} \tag{3.5}$$

Therefore, for a given $n$, $C_{u,n} \ge D_{u,n}$ if and only if $\sum_{j=3}^{n-1} D_{u,j} \ge D_{u,n}$. In fact, the inequality is strict for all sufficiently large $n$ as shown below in Theorem 3.1.

From the singularity analysis in Appendix B, we find $D_{u,n} = \frac{2^n}{4\sqrt{\pi n}}\{1 - \frac{d_u}{16n} + o(n^{-1})\}$ with constants $d_u = (4u+3)^2$ for $u$=0, 1, 2. Thus

$$\begin{aligned}D_{u,n}/D_{u,n-1} &= 2(1-\frac{1}{n})^{1/2}\left\{\frac{1+\frac{d_u}{16n}+o(n^{-1})}{1+\frac{d_u}{16(n-1)}+o(n^{-1})}\right\} = 2(1-\frac{1}{n})^{1/2}\{1+O(n^{-2})\} = 2(1-\frac{1}{2n}+o(n^{-1}))\{1+o(n^{-1})\} \\ &= 2-\frac{1}{n}+o(n^{-1}),\end{aligned}$$

so that for sufficiently large $n$,

$$D_{u,n}/D_{u,n-1} < 2 - \frac{1}{2n}, \text{ say,} \tag{3.6}$$

and

$$D_{u,n-1}/D_{u,n} = \frac{1}{2}\{1-\frac{1}{2n}+o(n^{-1})\}^{-1} = \frac{1}{2}\{1+\frac{1}{2n}+o(n^{-1})\} > \frac{1}{2}. \tag{3.7}$$

Theorem 3.1. *For u=0, 1, and 2, and n sufficiently large,*

$$\sum_{j=3}^{n-1} D_{u,j} > D_{u,n} . \tag{3.8}$$

Proof. Suppressing $u$ from the subscripts for this proof,

$$\sum_{j=3}^{n-1} D_j = D_{n-1}\left\{1+\frac{D_{n-2}}{D_{n-1}}+\frac{D_{n-3}}{D_{n-2}}\frac{D_{n-2}}{D_{n-1}}+\cdots+\frac{D_3\cdots D_{n-2}}{D_4\cdots D_{n-1}}\right\}$$

$$= D_{n-1}\sum_{j=1}^{n-3}\prod_{i=1}^{j-1}\frac{D_{n+i-j-1}}{D_{n+i-j}} > D_{n-1}\sum_{j=1}^{\lfloor n/2\rfloor}\prod_{i=1}^{j-1}\frac{D_{n+i-j-1}}{D_{n+i-j}}$$

since all $D_j$ are positive and $n-3>n/2$ for $n\ge 7$. Thus for each $i$, $j$ satisfying $1\le i<j\le n/2$ we have $n+i-j>n/2$, so by (3.7),

$$D_{n+i-j-1}/D_{n+i-j} > \frac{1}{2} \tag{3.9}$$

for sufficiently large $n$, whence

$$\sum_{j=3}^{n-1} D_j > D_{n-1}\sum_{j=1}^{\lfloor n/2\rfloor}\left(\frac{1}{2}\right)^{j-1} = D_{n-1}\{2-O(2^{-n/2})\} = D_{n-1}\{2+o(n^{-1})\} .$$

Thus by (3.6), $\sum_{j=3}^{n-1} D_j - D_n > D_{n-1}\{2+o(n^{-1})-2+\frac{1}{2n}\} = \frac{D_{n-1}}{2n}\{1+o(1)\} > \frac{D_{n-1}}{3n}$, say, for sufficiently large $n$, so for such $n$,

$$\sum_{j=3}^{n-1} D_j - D_{n-1} > 0 . \quad \square \tag{3.10}$$

With respect to smaller values of $n$, exact evaluation of $C_{u,n}$ and $D_{u,n}$ suggests the following values of $n$ qualify as "sufficiently large" for the argument of Lemma 3.1 to hold.

| Inequality | Eqn | Holds for $n$ at least | | |
|---|---|---|---|---|
| | | $u=0$ | $u=1$ | $u=2$ |
| $D_{u,n}/D_{u,n-1} < 2-\frac{1}{2n}$ | (3.6) | 17 | 36 | 53 |
| $D_{n+i-j-1}/D_{n+i-j} > \frac{1}{2}$ | (3.9) | 15 | 30 | 44 |
| $\sum_{j=3}^{n-1} D_{u,j} - D_{u,n-1} > \frac{D_{u,n-1}}{3n}$ | (3.10) | 11 | 27 | 41 |

Inequalities (3.6), (3.9), and (3.10) have been confirmed numerically for each value of $n$ from the tabled values up to $n$=500, for $u$=0, 1, and 2. As this extends well into the asymptotic regime, we conclude (3.8) holds at least for all $n \geq 17$, 36, and 53, respectively, and numerical evaluation confirms (3.8) also holds for $n$ as small as 8, 13, and 23, respectively..

This concludes the proof of Theorem 3.1. □

Remark 3.1. Of course, condition (1.12) may hold for smaller values of $n$ than shown above due to the $\Delta_{u,n}^{(3)}$ term in (1.14), which is positive for $u=1$ and 2. Indeed, under (1.10), we find (1.12) and (1.15) hold for $n$ as small as 7, 7, and 9 for $u=0$, 1, and 2, respectively, and strictly so except for $n=7$, $u=0$ where $C_{0,7} = 1\cdot 4 + 3\cdot 1 = 1\cdot 1 + 3\cdot 2 = D_{0,7}$ [and similarly for (1.11) where the equality is $a_1 + a_2 = 2\cdot 4 + 6\cdot 1 = 14 = 1\cdot 5 + 3\cdot 3 = b_1 + b_2$]. In all cases, though, $H_s(n) > H_{s-1}(n)$ for $n \geq 7$. This concludes the proof of Theorem 1.1. □

## 4. DISCUSSION

Table 4.1 summarizes the necessary and/or sufficient conditions which imply Theorem 1.1. In the table, we use $s$ in the subscript of probability measures for $s \geq 0$ and $u = -s$ in the subscript for $s \leq 0$. (The notation has been changed slightly from prior sections for clarity.) Figure 2 below diagrams the logical relations among the necessary and/or sufficient conditions.

[TABLE 4.1 HERE]

Table 4.1 Necessary and/or sufficient conditions for Theorem 1.1 to hold.

| Probability measures over $k=1,\ldots,\lfloor n_s/3\rfloor=\lfloor (n-s-1)/3\rfloor$ for $s\ge 0$. | |
|---|---|
| Single-factor measure<br>$P_{s,n}^{(1)}(k)\propto\binom{2k+s}{k}$ | Double-factor measure<br>$P_{s,n}^{(2)}(k)\propto\binom{2k+s}{k}\binom{n_s-2k}{k}$ |
| Necessary and sufficient condition<br>[NS1$s$] $E_{s,n}^{(1)}\left[\binom{n_s-1-2K}{K-1}\right]+\Delta_{s,n}^{(1)}\ \ge\ E_{s,n}^{(1)}\left[\binom{n_s-1-2K}{K}\right]-\Delta_{s,n}^{(1)}$ | Necessary and sufficient condition<br>[NS2$s$] $E_{s,n}^{(2)}\left[\left(\frac{K}{K+s+1}\right)\left(\frac{4K-\{n_s-(s+1)\}}{n_s-2K}\right)\right]\ge 0$ |
| Sufficient condition<br>----- | Sufficient condition<br>[S2$s$] $E_{s,n}^{(2)}\left[K\right]\ge\{n_s-(s+1)\}/4$ |
| Probability measures over $k=u+1,\ldots,\lfloor n_u/3\rfloor=\lfloor (n+u)/3\rfloor$ for $u=-s\ge 0$. | |
| Single-factor measure<br>$P_{u,n}^{(1)}(k)\propto\binom{2k-u-1}{k}$ | Double-factor measure<br>$P_{u,n}^{(2)}(k)\propto\binom{2k-u-1}{k}\binom{n+u-1-2k}{k-1}$ |
| Necessary and sufficient condition<br>[NS1$u$] $E_{u,n}^{(1)}\left[\binom{n_u-1-2K}{K}\right]+\Delta_{u,n}^{(1)}\ge E_{u,n}^{(1)}\left[\binom{n_u-1-2K}{K-1}\right]-\Delta_{u,n}^{(1)}$ | Necessary and sufficient condition<br>[NS2$u$] $E_{u,n}^{(2)}\left[(K-u)^{-1}\right]\ge 4/(n-2u)$ |
| Sufficient condition<br>[S1$u$] $E_{u,n}^{(1)}\left[\binom{n_u-1-2K}{K}\right]\ \ge\ E_{u,n}^{(1)}\left[\binom{n_u-1-2K}{K-1}\right]$ | Sufficient conditions<br>[S2$u$] $E_{u,n}^{(2)}\left[K\right]\le(n_u+u)/4$<br>[S2(1)$u$] $E_{u,n}^{(2)}\left[K^{-1}\right]\ge 4/n_u$ |
| Shifted probability measures over $k=1,\ldots,\lfloor n_u/3\rfloor-u=\lfloor (n-2u)/3\rfloor$ for $u=-s\ge 0$. | |
| Single-factor measure<br>$P_{u,n}^{(1)\prime}(k)\propto\binom{2k+u-1}{k-1}$ | Double-factor measure<br>$P_{u,n}^{(2)\prime}(k)\propto\binom{2k+u-1}{k-1}\binom{n-u-1-2k}{k+u-1}$ |
| Necessary and sufficient condition<br>[NS1$u$]′ $E_{u,n}^{(1)\prime}\left[\binom{n-u-1-2K}{K+u}\right]+\Delta_{u,n}^{(1)\prime}\ge E_{u,n}^{(1)\prime}\left[\binom{n-u-1-2K}{K+u-1}\right]-\Delta_{u,n}^{(1)\prime}$ | Necessary and sufficient condition<br>[NS2$u$]′ $E_{u,n}^{(2)\prime}\left[K^{-1}\right]\ge 4/(n-2u)$ |
| Sufficient condition<br>[S1$u$]′ $E_{u,n}^{(1)\prime}\left[\binom{n-u-1-2K}{K+u}\right]\ \ge\ E_{u,n}^{(1)\prime}\left[\binom{n-u-1-2K}{K+u-1}\right]$ | Sufficient conditions<br>[S2$u$]′ $E_{u,n}^{(2)\prime}\left[K\right]\le(n-2u)/4$<br>[S2(1)$u$]′ $E_{u,n}^{(2)\prime}\left[(K+u)^{-1}\right]\ge 4/n_u$ |

The definitions of the $\Delta$-offsets in the single-factor necessary and sufficient conditions are as follows.

$$\Delta_{s,n}^{(1)}=E_{s,n}^{(1)}\left[\left(\frac{s+1}{K+s+1}\right)\binom{n_s-1-2K}{K}\right],\ \Delta_{u,n}^{(1)}=E_{u,n}^{(1)}\left[\left(\frac{u}{K-u}\right)\binom{n_u-1-2K}{K}\right],\ \text{and}\ \Delta_{u,n}^{(1)\prime}=E_{u,n}^{(1)\prime}\left[\left(\frac{u}{K}\right)\binom{n-u-1-2K}{K+u}\right].$$

Figure 2. Logical relations amongst the conditions

$$s \geq 0$$

$$\begin{array}{ccc} [NS1s] & \Leftrightarrow & [NS2s] \\ & & \Uparrow \\ & & \mathbf{[S2s]} \end{array}$$

$$s = -u \leq 0$$

$$\begin{array}{ccc} [NS1u] & \Leftrightarrow & [NS2u] \\ \Uparrow & & \Uparrow \\ \vdots & & \mathbf{[S2u]} \\ \Uparrow & & \\ \mathbf{[S1u]} & \Leftrightarrow & [S2(1)u] \end{array}$$

Notes: Conditions in boldface are used to establish Theorem 1.1. An analogous lower panel holds for conditions [$NS1u$]′, [$S1u$]′, etc. with shifted domain $k = 1,...,\lfloor (n-2u)/3 \rfloor$.

Theorem 2.1 establishes sufficient condition $[S2s]$ which implies the necessary and sufficient condition $[NS2s]$ by Jensen's inequality (hence also the equivalent necessary and sufficient condition $[NS1s]$).

Theorem 3.1 establishes sufficient condition $[S1u]$ for $u = 0, 1, 2$ and *n* sufficiently large, and $[S1u]$ implies necessary and sufficient condition $[NS1u]$ *a fortiori*. Numerical evaluation confirms $[NS1u]$ directly for small *n* per Remark 3.1. Prior work (L25, Section 6) shows $[S2u]$ holds for $u > 2$ for all $n \geq 7$.

Although not needed for the proofs, we note that sufficient condition $[S2(1)u]$ is equivalent to $[S1u]$ by a simple change of measure argument via the identity $\binom{n+u-1-2k}{k} \Big/ \binom{n+u-1-2k}{k-1} = \frac{n+u-3k}{k} = \frac{n_u}{k} - 3$, for then [$S1u$]′ holds if and only if $E^{(2)}_{u,n}\left[\binom{n+u-1-2K}{K} \Big/ \binom{n+u-1-2K}{K-1}\right] \geq 1$ if and only if $E^{(2)}_{u,n}\left[\frac{n_u}{K} - 3\right] \geq 1$, which is $[S2(1)u]$. While another application of Jensen's inequality would yield a further sufficient condition implying [S2(1)*u*], namely, $E^{(2)'}_{u,n}[K] \leq n_u / 4$, almost no values of *n* satisfy that condition for any *u*.

**ACKNOWLEDGMENT**

The author would like to thank Cheng-Shiun Leu and Jacob Fink for helpful discussions.

**DISCLOSURE STATEMENT**

The author has no conflicts of interest to report. Gemini 3 Flash was used to look up known generating functions and Claude AI was used to cross-check my asymptotic calculations.

**APPENDIX A**

We derive formulas (2.6) and (2.7) using the transfer theorems of singularity analysis (Flajolet and Sedgewick, 2009) for generating functions of the form $G_{\alpha,\beta}(x) = F(x)(1-x/\beta)^{-\alpha}$ for positive constant $\beta$, real $\alpha \neq 0,-1,-2,\ldots$, and $F(x)$ analytic at $\beta$. For non-positive integer $\alpha$, the coefficients of $x^n$ in the polynomial generating function are zero for $n$ sufficiently large, so can be ignored in the asymptotic expansions. The basic transfer theorem states that the coefficient of $x^n$ in $(1-x/\beta)^{-\alpha}$, written $[x^n](1-x/\beta)^{-\alpha}$, is given by

$$[x^n](1-x/\beta)^{-\alpha} = \frac{n^{\alpha-1}\beta^{-n}}{\Gamma(\alpha)}\left\{1+\frac{\alpha(\alpha-1)}{2n}+\frac{\alpha(\alpha-1)(\alpha-2)(3\alpha-1)}{24n^2}+O(n^{-3})\right\}. \tag{A.1}$$

For $G_{\alpha,\beta}(x)$, expand $F(x)$ in a Taylor series around $\beta$:

$$\begin{aligned} F(x) &= F(\beta)+F'(\beta)(x-\beta)+½F''(\beta)(x-\beta)^2+O(x-\beta)^3 \\ &= F(\beta)-\beta F'(\beta)(1-x/\beta)+½\beta^2F''(\beta)(1-x/\beta)^2+O(1-x/\beta)^3 . \end{aligned}$$

Then

$$G_{\alpha,\beta}(x) = F(\beta)(1-x/\beta)^{-\alpha} - \beta F'(\beta)(1-x/\beta)^{-(\alpha-1)} + ½\beta^2F''(\beta)(1-x/\beta)^{-(\alpha-2)} + O(1-x/\beta)^{-(\alpha-3)}$$

and Theorem VI.1 of Flajolet and Sedgewick (2009) states that the asymptotic expansion of coefficients of $x^n$ from a linear combination of generating functions is the same linear combination of the individual expansions, so that

$$\begin{aligned} [x^n]F(x)(1-x/\beta)^{-\alpha} &= \frac{F(\beta)n^{\alpha-1}\beta^{-n}}{\Gamma(\alpha)}\left\{1+\frac{\alpha(\alpha-1)}{2n}+\frac{\alpha(\alpha-1)(\alpha-2)(3\alpha-1)}{24n^2}+O(n^{-3})\right\} \\ &\quad - \frac{\beta F'(\beta)n^{\alpha-2}\beta^{-n}}{\Gamma(\alpha-1)}\left\{1+\frac{(\alpha-1)(\alpha-2)}{2n}+O(n^{-2})\right\} \\ &\quad + \frac{½\beta^2F''(\beta)n^{\alpha-3}\beta^{-n}}{\Gamma(\alpha-2)}\left\{1+O(n^{-1})\right\}, \end{aligned}$$

whence

$$[x^n]F(x)(1-x/\beta)^{-\alpha} = \frac{n^{\alpha-1}\beta^{-n}}{\Gamma(\alpha)}\left\{ F(\beta)+c_1 n^{-1}+c_2 n^{-2}+O(n^{-3}) \right\} \tag{A.2}$$

where

$$c_1 = F(\beta)\frac{\alpha(\alpha-1)}{2} - \beta F'(\beta)\frac{\Gamma(\alpha)}{\Gamma(\alpha-1)} = (\alpha-1)\left\{\frac{\alpha}{2}F(\beta)-\beta F'(\beta)\right\} \tag{A.3}$$

and

$$\begin{aligned} c_2 &= F(\beta)\frac{\alpha(\alpha-1)(\alpha-2)(3\alpha-1)}{24} - \beta F'(\beta)\frac{(\alpha-1)(\alpha-2)}{2}\frac{\Gamma(\alpha)}{\Gamma(\alpha-1)} + ½\beta^2 F''(\beta)\frac{\Gamma(\alpha)}{\Gamma(\alpha-2)} \\ &= \frac{(\alpha-1)(\alpha-2)}{2}\left\{\frac{\alpha(3\alpha-1)}{12}F(\beta)-(\alpha-1)\beta F'(\beta)+\beta^2 F''(\beta)\right\}. \end{aligned} \tag{A.4}$$

Below we apply (A.2)–(A.4) to some standard combinatorial generating functions which we re-derive here for completeness. The reader may skip these derivations by jumping to "Transfer" below the proof of Lemma A.3.

Let $t \ge 0$ be a given integer and define $G^{(t)}(x) = x^t \sum_{k\ge 0}\binom{2k+t}{k}x^k$ so that the g.f. of the sequence $\left\{\binom{2k+t}{k} : k \ge 0\right\}$ is $x^{-t}G^{(t)}(x)$.

Lemma A.1. *(i)* $G^{(0)}(x) = (1-4x)^{-1/2}$.

*(ii)* $G^{(1)}(x) = ½\{(1-4x)^{-1/2}-1\}$.

*(iii)* $G^{(t)}(x) = \left\{\frac{1-(1-4x)^{1/2}}{2}\right\}^t (1-4x)^{-1/2}$ *for* $t > 1$.

Proof.

(i) $$(1-4x)^{-1/2} = \sum_{k\ge 0}\binom{-½}{k}(-4)^k x^k = \sum_{k\ge 0}\frac{(-1/2)(-3/2)\cdots\{-(2k-1)/2\}}{k!}(-4)^k x^k$$

$$= \sum_{k\ge 0}\frac{(1\cdot 3\cdots(2k-1)\}}{2^k k!}4^k\,\frac{(2\cdot 4\cdots(2k)\}}{2^k k!}x^k = \sum_{k\ge 0}\binom{2k}{k}x^k = G^{(0)}(x).$$

(ii) $\binom{2k+2}{k+1} = \binom{2k+1}{k+1}+\binom{2k+1}{k} = 2\binom{2k+1}{k}$, so

$$G^{(1)}(x) = x\sum_{k\ge 0}\binom{2k+1}{k}x^k = ½\sum_{k\ge 0}\binom{2k+2}{k+1}x^{k+1} = ½\sum_{k\ge 1}\binom{2k}{k}x^k = ½\{G^{(0)}(x)-1\} = ½\{(1-4x)^{-1/2}-1\}.$$

(iii) $G^{(t)}(x)$ satisfies a functional linear recursion in $t$, namely, $G^{(t+1)} = G^{(t)} - x\,G^{(t-1)}$, as follows.

$$G^{(t+1)}(x) = x^{t+1}\sum_{k\ge 0}\binom{2k+t+1}{k}x^k \ ,$$

$$xG^{(t-1)}(x) = x^t + x^t\sum_{k\ge 1}\binom{2k+t-1}{k}x^k = x^t + x^{t+1}\sum_{k\ge 0}\binom{2k+t+1}{k+1}x^k \ ,$$

so adding,

$$G^{(t+1)}(x) + xG^{(t-1)}(x) = x^t + x^{t+1}\sum_{k\ge 0}\left\{\binom{2k+t+1}{k+1}+\binom{2k+t+1}{k}\right\}x^k = x^t + x^{t+1}\sum_{k\ge 0}\binom{2k+t+2}{k+1}x^k$$

$$= x^t + x^t\sum_{k\ge 1}\binom{2k+t}{k}x^k = G^{(t)}(x) \, ,$$

which yields the recursion $G^{(t+1)} = G^{(t)} - xG^{(t-1)}$ whose characteristic equation $z^2 - z + x = 0$ has roots $r_1(x) = ½\{1-(1-4x)^{1/2}\}$ and $r_2(x) = ½\{1+(1-4x)^{1/2}\}$. Thus $G^{(t)} = C_1 r_1^t + C_2 r_2^t$ where $C_1$ and $C_2$ are functions that must satisfy the initial conditions $C_1 + C_2 = G^{(0)}$ and $C_1 r_1 + C_2 r_2 = G^{(1)}$. But $G^{(1)}(x) = ½\{(1-4x)^{-1/2} - 1\} = (1-4x)^{-1/2}\cdot ½\{1-(1-4x)^{1/2}\} = G^{(0)}(x) r_1(x)$, so solving for $C_1(x)$ and $C_2(x)$ we find $C_1(x) = G^{(0)}(x)$ and $C_2(x) = 0$. Thus $G^{(t)}(x) = r_1(x)^t(1-4x)^{-1/2}$ which is (iii). □

Lemma A.2. *The g.f. of* $\left\{\binom{n+\alpha-1}{\alpha-1} : n \ge 0\right\}$ *is* $(1-x)^{-\alpha}$.

Proof.

$$(1-x)^{-\alpha} = \sum_{n\ge 0}\binom{-\alpha}{n}(-1)^n x^n = \sum_{n\ge 0}\frac{(-\alpha)\cdots\{-\alpha-n+1\}}{n!}(-1)^n x^n = \sum_{n\ge 0}\frac{(n+\alpha-1)^{(n)}}{n!}x^n$$

$$= \sum_{n\ge 0}\binom{n+\alpha-1}{n}x^n = \sum_{n\ge 0}\binom{n+\alpha-1}{\alpha-1}x^n \, . \quad □$$

We can now obtain the g.f.'s of $A_n$ and $B_n$ defined at (2.3) and (2.4).

Lemma A.3.

*The g.f. of the sequence* $A_n = \left\{\sum_{k\ge 0}\binom{2k}{k}\binom{n-2k}{k} : n \ge 0\right\}$ *is*

$$G_A(x) = F_A(x)(1-2x)^{-1/2} \quad \textit{where} \quad F_A(x) = (1-x)^{-1/2}(1+x+2x^2)^{-1/2} \tag{A.5}$$

*and the g.f. of the sequence* $B_n = \left\{\sum_{k\ge 0}k\binom{2k}{k}\binom{n-2k}{k} : n \ge 0\right\}$ *is*

$$G_B(x) = F_B(x)(1-2x)^{-3/2} \quad \textit{where} \quad F_B(x) = 2x^3(1-x)^{-1/2}(1+x+2x^2)^{-3/2} \, . \tag{A.6}$$

Proof. $G_A(x) = \sum_{n\ge 0}\sum_{k\ge 0}\binom{2k}{k}\binom{n-2k}{k}x^n = \sum_{k\ge 0}\binom{2k}{k}\sum_{n\ge 3k}\binom{n-2k}{k}x^n = \sum_{k\ge 0}\binom{2k}{k}x^{3k}\sum_{n\ge 0}\binom{n+k}{k}x^n$

$$= \sum_{k\ge 0}\binom{2k}{k} x^{3k}(1-x)^{-(k+1)} \text{ by Lemma A.2 with } \alpha = k+1$$

$$= (1-x)^{-1}\sum_{k\ge 0}\binom{2k}{k}\left\{\frac{x^3}{1-x}\right\}^k = (1-x)^{-1}(1-4y)^{-1/2} \text{ by Lemma A.1(i),}$$

where $y = x^3/(1-x)$. But $1-4y = (1-x)^{-1}(1-x-4x^3) = (1-x)^{-1}(1+x+2x^2)(1-2x)$, so

$$G_A(x) = (1-x)^{-1}(1-x)^{1/2}(1+x+2x^2)^{-1/2}(1-2x)^{-1/2} = F_A(x)(1-2x)^{-1/2}.$$

For $G_B(x)$, we begin as above, leading to $G_B(x) = (1-x)^{-1}\sum_{k\ge 0} k\binom{2k}{k} y^k$, with $y = x^3/(1-x)$, so

$$G_B(x) = (1-x)^{-1} y \frac{d}{dy}\sum_{k\ge 0}\binom{2k}{k} y^k = (1-x)^{-1} y \frac{d}{dy}(1-4y)^{-1/2} = 2(1-x)^{-1} y\,(1-4y)^{-3/2}$$

$$= 2x^3(1-x)^{-2}(1-x)^{3/2}(1+x+2x^2)^{-3/2}(1-2x)^{-3/2} = F_B(x)(1-2x)^{-3/2}. \quad \square$$

Transfer. Now we apply (A.2)–(A.4) to $G_A(x)$ for the asymptotic expansion of the $A_n$ sequence.

$$F_A(x) = (1-x)^{-1/2}(1+x+2x^2)^{-1/2}$$

$$\log F_A(x) = -\tfrac{1}{2}\log(1-x) - \tfrac{1}{2}\log(1+x+2x^2)$$

$$\frac{d}{dx}\log F_A(x) = \tfrac{1}{2}(1-x)^{-1} - \tfrac{1}{2}(1+4x)(1+x+2x^2)^{-1}$$

$$\frac{d^2}{dx^2}\log F_A(x) = \tfrac{1}{2}(1-x)^{-2} - \tfrac{1}{2}\{4(1+x+2x^2)^{-1} - (1+4x)^2(1+x+2x^2)^{-2}\}$$

$$F_A'(x) = F_A(x)\frac{d}{dx}\log F_A(x)$$

$$F_A''(x) = F_A(x)[\frac{d^2}{dx^2}\log F_A(x) + \{\frac{d}{dx}\log F_A(x)\}^2]$$

At $\beta = \tfrac{1}{2}$, $F_A(\beta) = 1$, $\frac{d}{dx}\log F_A(\beta) = \frac{1}{2}\cdot 2 - \frac{1}{2}\cdot 3\cdot 2^{-1} = \frac{1}{4}$, $F_A'(\beta) = \frac{1}{4}$, and

$$\frac{d^2}{dx^2}\log F_A(\beta) = \frac{1}{2}\cdot 4 - \frac{1}{2}(4\cdot\frac{1}{2} - 9\cdot\frac{1}{4}) = \frac{17}{8}, \text{ so that } F_A''(\beta) = \frac{17}{8} + \frac{1}{16} = \frac{35}{16}.$$

Thus (A.3) at $\alpha = \tfrac{1}{2}$ is

$$c_1 = (\alpha-1)\{\frac{\alpha}{2}F(\beta) - \beta F'(\beta)\} = -\frac{1}{2}(\frac{1}{4}\cdot 1 - \frac{1}{2}\cdot\frac{1}{4}) = \frac{-1}{16} \quad \text{and}$$

$$c_2 = \frac{(\alpha-1)(\alpha-2)}{2}\{\frac{\alpha(3\alpha-1)}{12}F_A(\beta) - (\alpha-1)\beta F_A'(\beta) + \beta^2 F_A''(\beta)\}$$

$$= \frac{1}{2}(-\frac{1}{2})(-\frac{3}{2})\{\frac{1}{12}\cdot\frac{1}{2}\cdot\frac{1}{2}\cdot 1 \;-\; (-\frac{1}{2})\frac{1}{2}\cdot\frac{1}{4} \;+\; \frac{1}{4}\cdot\frac{35}{16} = \frac{121}{512}.$$

Therefore,

$$A_n = [x^n]F_A(x)(1-2x)^{-1/2} = \frac{2^n n^{-1/2}}{\sqrt{\pi}}\{1 - \frac{1}{16}n^{-1} + \frac{121}{512}n^{-2} + O(n^{-3})\} \quad \text{which is (2.6).}$$

Next we apply (A.2)–(A.4) to $G_B(x)$ for the asymptotic expansion of the $B_n$ sequence.

$$F_B(x) = 2x^3(1-x)^{-1/2}(1+x+2x^2)^{-3/2}$$

$$\log F_B(x) = \log 2 + 3\log x - ½\log(1-x) - \frac{3}{2}\log(1+x+2x^2)$$

$$\frac{d}{dx}\log F_B(x) = 3x^{-1} + ½(1-x)^{-1} - \frac{3}{2}(1+4x)(1+x+2x^2)^{-1}$$

$$\frac{d^2}{dx^2}\log F_B(x) = -3x^{-2} + ½(1-x)^{-2} - \frac{3}{2}\{4(1+x+2x^2)^{-1} - (1+4x)^2(1+x+2x^2)^{-2}\}$$

$$F_B'(x) = F_B(x)\frac{d}{dx}\log F_B(x)$$

$$F_B''(x) = F_B(x)[\frac{d^2}{dx^2}\log F_B(x) + \{\frac{d}{dx}\log F_B(x)\}^2]$$

At $\beta = ½$, we have $F_B(\beta) = 2\cdot 2^{-3}\cdot 2^{1/2}\cdot 2^{-3/2} = \frac{1}{8}$,

$\frac{d}{dx}\log F_B(\beta) = 3\cdot 2 + ½\cdot 2 - \frac{3}{2}\cdot 3\cdot 2^{-1} = \frac{19}{4}$, and $\frac{d^2}{dx^2}\log F_B(\beta) = -12 + 2 - \frac{3}{2}(4\cdot 2^{-1} - 9\cdot 2^{-2}) = \frac{-77}{8}$,

so that $F_B'(\beta) = \frac{1}{8}\cdot\frac{19}{4} = \frac{19}{32}$ and $F_B''(\beta) = \frac{1}{8}\{\frac{-77}{8} + (\frac{19}{4})^2\} = \frac{207}{128}$.

Thus with $\alpha = 3/2$, $c_1 = (\alpha-1)\{\frac{\alpha}{2}F(\beta) - \beta F'(\beta)\} = (\frac{3}{2}-1)(\frac{3}{4}\cdot\frac{1}{8} - \frac{1}{2}\cdot\frac{19}{32}) = \frac{-13}{128}$ and

$$c_2 = \frac{(\alpha-1)(\alpha-2)}{2}\{\frac{\alpha(3\alpha-1)}{12}F_B(\beta) - (\alpha-1)\beta F_B'(\beta) + \beta^2 F_B''(\beta)\}$$

$$= \frac{1}{2}\cdot(\frac{3}{2}-1)(\frac{3}{2}-2)\{\frac{1}{12}\cdot\frac{3}{2}\cdot(\frac{9}{2}-1)\cdot\frac{1}{8} - (\frac{3}{2}-\frac{1}{2})\cdot\frac{1}{2}\cdot\frac{19}{32} + \frac{1}{4}\cdot\frac{207}{128} = \frac{-159}{4096}.$$

so, $B_n = [x^n]F_B(x)(1-2x)^{-3/2} = \frac{(1/8)2^n n^{1/2}}{(1/2)\sqrt{\pi}}\{1 - \frac{13}{16}n^{-1} - \frac{159}{512}n^{-2} + O(n^{-3})\}$, which is (2.7).

For the asymptotic expansion of $B_n / A_n$, we have

$$\frac{B_n}{A_n} = \frac{2^n n^{1/2}}{4\sqrt{\pi}}\left\{1 - \frac{13}{16n} - \frac{159}{512n^{-2}} + O(n^{-3})\right\} \Big/ \frac{2^n n^{-1/2}}{\sqrt{\pi}}\left\{1 - \frac{1}{16n} + \frac{121}{512n^{-2}} + O(n^{-3})\right\}$$

$$= \frac{n}{4}\left\{1 - \frac{13}{16n} - \frac{159}{512n^{-2}} + O(n^{-2})\right\}\left\{1 + \left(\frac{1}{16n} - \frac{121}{512n^{-2}}\right) + \left(\frac{1}{16n} - \frac{121}{512n^{-2}}\right)^2 + O(n^{-3})\right\}$$

$$= \frac{n}{4}\left\{1 - \frac{3}{4}n^{-1} - \left(\frac{159}{512} + \frac{13}{256} + \frac{121}{512} - \frac{1}{256}\right)n^{-2} + O(n^{-3})\right\} = \frac{n}{4}\left\{1 - \frac{3}{4}n^{-1} - \frac{19}{32}n^{-2} + O(n^{-3})\right\}$$

$$= \frac{n}{4} - \frac{3}{16} - \frac{19}{128n} + O(n^{-2}) .$$

**APPENDIX B**

We derive the generating functions of the sequences $C_{u,n} = \sum_{k\ge1}\binom{2k+u-1}{k-1}\binom{n-u-1-2k}{k+u}$ and $D_{u,n} = \sum_{k\ge1}\binom{2k+u-1}{k-1}\binom{n-u-1-2k}{k+u-1}$, starting with two more standard results.

<u>Lemma B.1</u>. *For* $u \ge 0$, *the g.f. of* $\left\{\binom{2k+u-1}{k-1} : k \ge 1\right\}$ *equals* $x\left\{\frac{1-(1-4x)^{1/2}}{2x}\right\}^{u+1}(1-4x)^{-1/2}$.

<u>Proof</u>. $\sum_{k\ge1}\binom{2k+u-1}{k-1}x^k = x\sum_{k\ge0}\binom{2k+u+1}{k}x^k = x \cdot x^{-(u+1)}G^{(u+1)}(x)$ with the function $G^{(u+1)}(x)$ of Lemma A.1(iii) with $t = u+1$. Thus $\sum_{k\ge1}\binom{2k+u-1}{k-1}x^k = x\left\{\frac{1-(1-4x)^{1/2}}{2x}\right\}^{u+1}(1-4x)^{-1/2}$. □

<u>Lemma B.2</u>. *For* $k \ge 1$ *and* $u \ge 0$, *the g.f. of* $\left\{\binom{n-u-1-2k}{k+u-1} : n \ge 2u+3k\right\}$ *equals* $\left(\frac{x^2}{1-x}\right)^u\left(\frac{x^3}{1-x}\right)^k$.

<u>Proof</u>. Since $\binom{n-u-1-2k}{k+u-1} = \binom{n-2u-3k+k+u-1}{k+u-1}$ we have from Lemma A.2 with $\alpha = k+u$

$$\sum_{n\ge 2u+3k}\binom{n-u-1-2k}{k+u-1}x^n = x^{2u+3k}\sum_{n\ge0}\binom{n+\alpha-1}{\alpha-1}x^n = x^{2u+3k}(1-x)^{-(k+u)} . \quad \square$$

<u>Lemma B.3</u>. *The g.f. of* $\{D_{u,n} : n \ge 3\}$ *equals* $G_u^{(D)}(x) = 2^{-(u+1)}x^{-u}\sum_{i=0}^{u+1}(-1)^i\binom{u+1}{i}(1-\frac{4x^3}{1-x})^{(i-1)/2}$.

<u>Proof</u>.

$$G_u^{(D)}(x)=\sum_{n\ge 3}\sum_{k\ge 1}\binom{2k+u-1}{k-1}\binom{n-u-1-2k}{k+u-1}x^n=\sum_{k\ge 1}\binom{2k+u-1}{k-1}\sum_{n\ge 2u+3k}\binom{n-u-1-2k}{k+u-1}x^n$$

$$=\left(\frac{x^2}{1-x}\right)^u\sum_{k\ge 1}\binom{2k+u-1}{k-1}y^k \text{ by Lemma B.2, where } y=\frac{x^3}{1-x}$$

$$=\left(\frac{x^2}{1-x}\right)^u y\left\{\frac{1-(1-4y)^{1/2}}{2y}\right\}^{u+1}(1-4y)^{-1/2} \text{ by Lemma B.1}$$

$$=2^{-(u+1)}x^{2u}(1-x)^{-u}x^{-3u}(1-x)^u\{1-(1-4y)^{1/2}\}^{u+1}(1-4y)^{-1/2}=2^{-(u+1)}x^{-u}\{1-(1-4y)^{1/2}\}^{u+1}(1-4y)^{-1/2}$$

$$=2^{-(u+1)}x^{-u}\sum_{i=0}^{u+1}(-1)^i\binom{u+1}{i}(1-4y)^{(i-1)/2}. \quad \square$$

It follows that

$$G_u^{(D)}(x)=\sum_{i=0}^{u+1}F_u^{(i)}(x)(1-2x)^{-(1-i)/2}$$

where (B.1)

$$F_u^{(i)}(x)=(-1)^i 2^{-(u+1)}\binom{u+1}{i}x^{-u}(1-x)^{(1-i)/2}(1+x+2x^2)^{(i-1)/2}.$$

Transfer. For the argument of Section 3, we only need terms up to order $n^{-1}$ in the factor multiplying the leading term $2^n/\sqrt{\pi n}$ in the asymptotic expansion of $D_{u,n}$, which means that we may ignore values of $i>2$ in (B.1), because such values of *i* transfer only terms of smaller order. We can also ignore $i=1$ (or any odd value of *i*) because such *i* only contribute integer powers $p\ge 0$ of $(1-2x)$ in $G_u^{(D)}(x)$ which transfer to $[x^n](1-2x)^p=0$ for large *n*. So instead of $G_u^{(D)}(x)$ we may consider the simpler g.f.

$$F_u^{(0)}(x)(1-2x)^{-1/2}+F_u^{(2)}(x)(1-2x)^{1/2}$$

which transfers to

$$D_{u,n}=\frac{2^n n^{-1/2}}{\Gamma(½)}\left[F_u^{(0)}(½)-½\{¼F_u^{(0)}(½)-½F_u^{(0)\prime}(½)\}n^{-1}+O(n^{-2})\right]+\frac{2^n n^{-3/2}}{\Gamma(-½)}\{F_u^{(2)}(½)+O(n^{-1})\}$$

$$=\frac{2^n}{\sqrt{\pi n}}\left[F_u^{(0)}(½)-½\{¼F_u^{(0)}(½)-½F_u^{(0)\prime}(½)+F_u^{(2)}(½)\}n^{-1}+O(n^{-2})\right] \quad \text{(B.2)}$$

since $\Gamma(½)=-½\Gamma(-½)$, so that $\Gamma(-½)=-2\Gamma(½)=-2\sqrt{\pi}$.

Now we evaluate the quantities in (B.2). From (B.1),

$$\log F_u^{(i)}(x)=const-u\log x+½(1-i)\log(1-x)+½(i-1)\log(1+x+2x^2)$$

$$\frac{d}{dx}\log F_u^{(i)}(x)=-ux^{-1}+½(i-1)\{(1-x)^{-1}+(1+4x)(1+x+2x^2)^{-1}\}$$

$$(-1)^i F_u^{(i)}(½) = 2^{-(u+1)}\binom{u+1}{i}(½)^{-u}(½)^{(1-i)/2}(½)^{(1-i)/2} = (½)^{u+1}\binom{u+1}{i}(½)^{1-i-u} = (½)^{2-i}\binom{u+1}{i}.$$

So for $i=0$, $F_u^{(0)}(½) = ¼$ for any $u$, and for $i=2$, $F_0^{(2)}(½) = 0$, $F_1^{(2)}(½) = 1$, and $F_2^{(2)}(½) = 3$, which we can write as $F_u^{(2)}(½) = ½u(u+1)$ for $u=0$, 1, and 2.

$$\frac{d}{dx}\log F_u^{(i)}(½) = -2u + ½(i-1)(2+3/2) = -\frac{8u+7}{4} \text{ for } i=0\text{, whence}$$

$$F_u^{(0)\prime}(½) = F_u^{(0)}(½)\frac{d}{dx}\log F_u^{(0)}(½) = -\frac{8u+7}{16}.$$

Then from (B.2),

$$D_{u,n} = \frac{2^n}{\sqrt{\pi n}}[\frac{1}{4} - \frac{1}{2}\{\frac{1}{4}\cdot\frac{1}{4} + \frac{8u+7}{32} + \frac{u(u+1)}{2}\}n^{-1} + O(n^{-2})]$$

$$= \frac{2^n}{4\sqrt{\pi n}}[1 - \{\frac{1}{8} + \frac{8u+7}{16} + u(u+1)\}n^{-1} + O(n^{-2})]$$

and since $16\{\frac{1}{8} + \frac{8u+7}{16} + u(u+1)\} = 2 + 8u + 7 + 16u(u+1) = 16u^2 + 24u + 9 = (4u+3)^2$,

$$D_{u,n} = \frac{2^n}{4\sqrt{\pi n}}\{1 - \frac{(4u+3)^2}{16n} + O(n^{-2})\}\text{, as asserted in Section 3.}$$

Finally, $G_u^{(C)}(x)$ and $G_u^{(D)}(x)$ are related as follows.

Lemma B.4. *The g.f. of* $\{C_{u,n} : n \geq 3\}$ *is* $G_u^{(C)}(x) = \frac{x}{1-x}G_u^{(D)}(x)$.

Proof. Since $\binom{n-u-1-2k}{k+u} = \binom{n-2u-3k-1+k+u}{k+u}$, we proceed as in Lemma B.2 with $\alpha = k+u+1$ to find the g.f. of $\left\{\binom{n-u-1-2k}{k+u} : n \geq 2u+3k+1\right\}$ equals $\frac{x^{2u+3k+1}}{(1-x)^{-(k+u+1)}} = \frac{x}{1-x}\left(\frac{x^2}{1-x}\right)^u\left(\frac{x^3}{1-x}\right)^k$. The rest of the derivation proceeds exactly as in Lemma B.3 with the extra factor of $\frac{x}{1-x}$. □

**APPENDIX C.**

The reader may find the following glossary of notation used in this paper helpful. Symbols are listed roughly in order of appearance.

| Symbol | Meaning |
|---|---|
| $x^{(n)}$ | An $n$-bit string in $\{0,1\}^n$. |
| $S(x^{(n)})$ | The net score in bit string $x^{(n)}$, equal to Alice's point total minus Bob's point total. |
| $H_s(n)$ | Number of heady-$s$ $n$-strings. |
| $T_s(n)$ | Number of taily-$s$ $n$-strings. |
| $n_s$ | $n-s-1$ |
| $n_u$ | $n+u$ |
| $\kappa$ | The largest support point for a given discrete distribution under discussion, typically $\kappa_s = \lfloor n_s/3 \rfloor$ for $s \ge 0$ or $\kappa_u = \lfloor n_u/3 \rfloor$ for $u = -s \ge 0$ in Section 1, or shifted $\kappa = \lfloor n_u/3 \rfloor - u = \lfloor (n-2u)/3 \rfloor$ in Section 3. |
| $K$ | A discrete random variable with support $1,\ldots,\kappa$ and probability function given in equation (1.5), (1.7), (1.13), or (1.16). |
| $\Delta_{s,n}^{(1)}$ | See equation (1.6). |
| $f(x)$ | Function whose expectation $E_{s,n}^{(2)}[f(K)]$ is to be shown positive. See just below equation (1.7). |
| $E_{s,n}^{(2)}[K]$ | Expected value of $K$ with lower bound given in equation (1.8) |
| $a_n, b_n$ | Abbreviations defined in equation (1.9). |
| $\Delta_{u,n}^{(3)}$ | See equation (1.14). |
| $E_{u,n}^{(4)}[K^{-1}]$ | Expected value of $K^{-1}$ with lower bound given in equation (1.15) |
| $A_n, B_n$ | Sums defined in equations (2.3) and (2.4). |
| $C_{u,n}, D_{u,n}$ | Sums defined in equations (3.3) and (3.4). |
| $G_u^{(C)}(x), G_u^{(D)}(x)$ | Generating functions for $C_{u,n}, D_{u,n}$, respectively. |
| $G_{\alpha,\beta}(x)$ | Generating function of the form $F(x)(1-x/\beta)^{-\alpha}$ with $F(x)$ analytic at $\beta$. |
| $[x^n]G(x)$ | The coefficient of $x^n$ in an asymptotic expansion of $G(x)$. |
| $G^{(t)}(x)$ | See Lemma A.1(iii). |
| $G_A(x), G_B(x)$ | Generating functinos for $A_n$ and $B_n$; see equations (A.5) and (A.6). |
| $F_A(x), F_B(x)$ | Functions analytic at $x=½$ used in $G_A(x), G_B(x)$, respectively. |
| $F_u^{(i)}(x)$ | Functions analytic at $x=½$ used in $G_u^{(D)}(x)$; see equation (B.1). |

Figure 1. Number of heady-$s$ (left) and taily-$s$ (right) 100-bit strings with score $S(x^{(100)}) = s$ (center)

| Heady binary 100-sequences | Score (A-B) | Taily binary 100-sequences |
|---:|---:|:---|
| 1 | 99 | 0 |
| 1 | 98 | 0 |
| 1 | 97 | 1 |
| 99 | 96 | 2 |
| 195 | 95 | 3 |
| 289 | 94 | 101 |
| 5,037 | 93 | 293 |
| 14,151 | 92 | 576 |
| 27,349 | 91 | 5,507 |
| 187,235 | 90 | 19,263 |
| 618,559 | 89 | 45,651 |
| 1,429,157 | 88 | 226,533 |
| 5,894,416 | 87 | 819,981 |
| 19,787,954 | 86 | 2,179,827 |
| 50,776,450 | 85 | 7,786,074 |
| 162,766,422 | 84 | 26,498,036 |
| 512,790,375 | 83 | 74,118,078 |
| 1,373,066,419 | 82 | 226,320,014 |
| 3,902,881,691 | 81 | 703,460,883 |
| 11,328,246,625 | 80 | 1,971,191,559 |
| 30,261,895,052 | 79 | 5,556,910,449 |
| 80,721,529,654 | 78 | 15,912,811,707 |
| 218,127,436,190 | 77 | 43,378,674,030 |
| 566,857,733,530 | 76 | 116,191,998,480 |
| 1,447,505,233,470 | 75 | 311,849,764,898 |
| 3,696,328,798,194 | 74 | 816,801,626,970 |
| 9,251,451,151,354 | 73 | 2,097,437,382,770 |
| 22,741,870,195,614 | 72 | 5,346,861,177,274 |
| 55,470,228,378,916 | 71 | 13,422,782,937,198 |
| 133,481,222,574,320 | 70 | 33,130,830,072,498 |
| 316,326,794,037,664 | 69 | 80,861,334,894,864 |
| 741,513,302,510,336 | 68 | 194,871,056,788,872 |
| 1,717,937,380,190,467 | 67 | 462,978,196,766,592 |
| 3,928,313,184,895,683 | 66 | 1,086,613,692,236,864 |
| 8,879,204,353,637,955 | 65 | 2,520,045,134,185,059 |
| 19,844,756,741,510,377 | 64 | 5,771,011,239,304,847 |
| 43,830,503,264,231,408 | 63 | 13,058,058,616,794,281 |
| 95,709,809,952,804,746 | 62 | 29,205,964,719,715,611 |
| 206,699,041,895,744,274 | 61 | 64,559,332,690,925,362 |
| 441,437,992,463,137,750 | 60 | 141,066,503,440,759,184 |
| 932,406,661,862,926,202 | 59 | 304,780,048,669,716,366 |
| 1,948,223,648,934,796,070 | 58 | 651,124,632,918,764,918 |
| 4,027,026,235,562,231,934 | 57 | 1,375,608,703,526,484,454 |
| 8,235,079,725,839,081,226 | 56 | 2,874,397,551,913,079,838 |
| 16,662,449,915,334,882,500 | 55 | 5,940,937,571,676,447,994 |
| 33,359,797,923,255,942,792 | 54 | 12,146,396,964,707,263,014 |
| 66,090,648,536,684,566,264 | 53 | 24,567,617,264,675,707,752 |
| 129,573,190,174,119,127,560 | 52 | 49,162,187,325,615,348,800 |
| 251,401,875,455,624,459,360 | 51 | 97,336,013,395,670,438,920 |
| 482,739,672,079,389,412,720 | 50 | 190,683,311,492,296,499,928 |
| 917,406,246,357,600,418,384 | 49 | 369,630,410,641,457,049,040 |
| 1,725,539,414,955,267,709,104 | 48 | 709,010,084,039,599,076,930 |
| 3,212,242,501,400,858,796,294 | 47 | 1,345,793,510,709,689,751,504 |
| 5,918,553,039,112,830,141,286 | 46 | 2,527,869,867,503,690,485,208 |
| 10,793,136,521,972,706,722,670 | 45 | 4,698,795,532,762,424,027,838 |
| 19,480,610,383,310,928,331,274 | 44 | 8,643,234,729,740,201,085,144 |
| 34,799,686,893,657,506,500,386 | 43 | 15,733,439,888,793,661,642,162 |
| 61,525,986,141,028,939,916,662 | 42 | 28,341,662,496,353,209,734,362 |
| 107,656,945,315,318,334,467,838 | 41 | 50,521,348,381,888,331,025,654 |
| 186,428,679,987,217,690,076,330 | 40 | 89,117,594,662,003,382,515,574 |
| 319,488,264,449,642,933,840,772 | 39 | 155,553,202,188,479,518,243,642 |
| 541,813,929,286,618,396,569,288 | 38 | 268,662,225,711,426,945,969,366 |
| 909,237,326,205,234,478,586,456 | 37 | 459,122,673,761,791,553,768,296 |
| 1,509,774,769,659,649,642,860,264 | 36 | 776,291,439,257,715,916,131,632 |
| 2,480,436,089,063,566,646,870,075 | 35 | 1,298,587,722,933,143,344,723,944 |
| 4,031,764,508,388,628,498,702,987 | 34 | 2,149,026,620,382,547,845,305,976 |
| 6,483,072,068,528,203,471,494,635 | 33 | 3,518,089,453,773,982,088,526,123 |
| 10,312,155,810,944,510,786,870,689 | 32 | 5,696,840,514,149,148,695,978,709 |
| 16,224,162,721,434,220,446,844,346 | 31 | 9,124,052,020,608,841,075,295,073 |
| 25,245,122,117,635,905,908,765,700 | 30 | 14,452,040,734,702,910,116,074,379 |
| 38,846,324,380,401,986,244,035,172 | 29 | 22,636,865,686,409,264,254,501,044 |
| 59,105,957,569,936,338,826,494,284 | 28 | 35,059,337,758,629,875,579,434,792 |
| 88,913,890,944,778,375,786,980,000 | 27 | 53,683,715,958,521,326,910,638,140 |
| 132,223,796,698,102,916,774,118,520 | 26 | 81,260,681,830,406,711,870,455,500 |
| 194,353,435,870,671,993,345,471,720 | 25 | 121,579,759,632,223,347,952,837,784 |
| 282,328,383,467,004,192,113,164,376 | 24 | 179,773,262,168,760,388,121,593,180 |
| 405,256,269,880,074,384,672,003,316 | 23 | 262,668,519,032,807,171,289,039,000 |
| 574,707,494,182,077,295,504,958,116 | 22 | 379,177,047,617,570,252,962,464,376 |
| 805,064,399,223,535,006,392,276,564 | 21 | 540,698,081,510,384,632,489,343,252 |
| 1,113,784,725,057,420,390,898,056,060 | 20 | 761,499,440,351,804,320,793,881,344 |
| 1,521,508,176,952,555,201,847,739,952 | 19 | 1,059,021,636,128,694,393,487,551,116 |
| 2,051,919,492,645,521,577,708,196,520 | 18 | 1,454,032,686,421,353,870,369,484,908 |
| 2,731,270,804,631,809,509,913,019,704 | 17 | 1,970,543,671,894,290,694,519,304,904 |
| 3,587,464,585,726,487,879,487,583,496 | 16 | 2,635,382,073,167,034,699,078,115,998 |
| 4,648,610,773,623,124,070,855,118,850 | 15 | 3,477,315,830,887,864,427,465,039,880 |
| 5,941,002,356,696,473,996,374,566,930 | 14 | 4,525,631,065,922,960,209,340,762,784 |
| 7,486,506,153,488,489,119,880,885,962 | 13 | 5,808,095,649,989,789,677,027,641,402 |
| 9,299,440,672,697,601,144,819,692,958 | 12 | 7,348,293,443,754,189,186,674,179,328 |
| 11,383,107,942,015,224,041,565,165,202 | 11 | 9,162,391,685,539,485,398,620,937,526 |
| 13,726,253,228,510,698,308,395,923,150 | 10 | 11,255,504,406,878,769,546,771,063,422 |
| 16,299,832,280,329,031,874,932,110,246 | 9 | 13,617,930,198,589,445,007,542,167,118 |
| 19,054,551,587,812,441,861,863,521,410 | 8 | 16,221,659,304,967,595,839,010,062,614 |
| 21,919,691,260,889,173,237,712,093,132 | 7 | 19,017,643,089,285,046,167,418,950,642 |
| 24,803,700,335,456,143,157,165,587,648 | 6 | 21,934,374,352,353,265,943,438,666,574 |
| 27,596,952,942,926,559,525,282,640,400 | 5 | 24,878,314,932,105,716,911,253,056,032 |
| 30,176,862,838,132,381,261,630,621,456 | 4 | 27,736,606,791,569,622,657,663,919,272 |
| 32,415,280,128,814,430,218,567,176,237 | 3 | 30,382,303,949,012,456,371,571,121,296 |
| 34,187,762,756,557,370,536,547,224,845 | 2 | 32,682,070,742,384,229,472,429,795,440 |
| 35,383,969,927,716,020,289,721,407,565 | 1 | 34,505,932,951,652,752,723,451,990,765 |
| 35,918,123,390,308,556,777,324,094,535 | 0 | 35,738,289,179,539,587,978,601,128,017 |
| 35,738,289,179,539,587,978,601,128,016 | −1 | 36,289,054,571,966,343,162,982,300,935 |
| 34,833,204,535,289,265,105,340,234,646 | −2 | 36,103,588,981,137,449,970,153,197,685 |
| 33,235,548,505,804,567,858,639,027,726 | −3 | 35,170,022,110,956,073,854,868,953,006 |
| 31,020,932,923,519,029,285,168,162,314 | −4 | 33,522,771,335,401,685,331,775,353,312 |
| 28,302,434,323,882,843,601,880,457,390 | −5 | 31,241,458,411,482,817,355,259,331,986 |
| 25,221,118,061,906,865,004,014,726,978 | −6 | 28,445,027,250,311,652,879,062,089,226 |
| 21,933,615,377,470,849,654,896,091,626 | −7 | 25,281,558,415,990,606,047,694,533,378 |
| 18,598,285,852,355,974,105,266,868,270 | −8 | 21,914,945,916,220,863,427,935,835,370 |
| 15,361,732,181,264,397,834,881,213,580 | −9 | 18,510,116,975,158,979,313,831,161,150 |
| 12,347,374,286,428,861,761,808,594,680 | −10 | 15,218,724,406,163,701,899,244,316,034 |
| 9,647,437,816,348,023,834,497,638,792 | −11 | 12,167,161,011,999,147,785,839,718,360 |
| 7,319,133,224,355,895,919,224,903,672 | −12 | 9,448,341,298,991,214,573,277,202,320 |
| 5,385,111,545,863,694,323,971,291,332 | −13 | 7,118,044,488,587,829,279,366,489,912 |
| 3,837,621,183,567,044,339,976,409,908 | −14 | 5,195,847,477,262,134,592,671,812,104 |
| 2,645,286,105,208,860,741,252,873,300 | −15 | 3,669,955,194,089,216,571,814,500,948 |
| 1,761,168,112,566,425,050,764,449,724 | −16 | 2,504,702,705,986,271,719,207,546,722 |
| 1,130,791,062,803,617,982,945,761,530 | −17 | 1,649,252,372,824,461,186,956,668,956 |
| 699,055,766,591,387,938,203,857,618 | −18 | 1,046,062,640,544,348,619,825,232,108 |
| 415,374,254,893,943,032,132,479,402 | −19 | 638,012,355,002,176,277,925,740,826 |
| 236,792,921,340,816,281,388,421,950 | −20 | 373,522,369,646,234,311,464,004,248 |
| 129,256,206,279,716,060,143,200,654 | −21 | 209,501,148,582,178,815,512,188,278 |
| 67,419,373,033,122,688,687,220,298 | −22 | 112,344,456,507,943,889,284,392,654 |
| 33,527,852,909,460,811,502,507,490 | −23 | 57,473,569,589,999,692,124,077,770 |
| 15,859,370,431,404,614,206,932,150 | −24 | 27,985,257,118,382,331,751,657,530 |
| 7,117,447,315,052,564,382,236,940 | −25 | 12,937,777,834,840,784,434,096,614 |
| 3,022,309,747,093,756,328,885,256 | −26 | 5,663,800,137,183,526,859,111,850 |
| 1,210,759,182,375,433,546,766,456 | −27 | 2,341,194,829,734,025,617,810,920 |
| 456,147,934,301,281,601,339,016 | −28 | 910,995,498,887,337,507,714,336 |
| 161,061,011,277,911,716,880,244 | −29 | 332,584,709,098,454,292,283,272 |
| 53,098,699,373,636,179,546,500 | −30 | 113,508,471,704,073,011,209,560 |
| 16,277,955,447,254,971,270,500 | −31 | 36,072,868,597,533,184,783,716 |
| 4,619,188,462,486,596,073,068 | −32 | 10,628,560,846,452,378,469,353 |
| 1,207,221,169,430,295,437,187 | −33 | 2,889,491,008,427,825,923,212 |
| 288,936,209,588,618,581,563 | −34 | 720,924,306,843,802,902,624 |
| 62,925,282,662,073,940,815 | −35 | 164,079,079,487,415,243,207 |
| 12,378,589,533,995,369,741 | −36 | 33,831,810,427,149,973,672 |
| 2,181,004,953,813,355,747 | −37 | 6,270,069,018,350,771,361 |
| 340,775,510,811,089,725 | −38 | 1,034,903,899,874,324,725 |
| 46,664,101,798,333,425 | −39 | 150,486,844,721,310,125 |
| 5,520,966,790,041,675 | −40 | 19,029,607,923,745,335 |
| 554,571,755,753,260 | −41 | 2,059,809,483,661,475 |
| 46,264,547,376,810 | −42 | 187,136,060,770,925 |
| 3,115,670,976,810 | −43 | 13,917,645,773,010 |
| 163,097,199,790 | −44 | 820,059,359,420 |
| 6,297,857,125 | −45 | 36,622,078,710 |
| 166,198,305 | −46 | 1,164,057,510 |
| 2,658,985 | −47 | 23,958,305 |
| 20,875 | −48 | 273,175 |
| 50 | −49 | 1,275 |
| 0 | −50 | 1 |